\documentclass{ws-jktr}
\usepackage{tikz-cd}
\usepackage{mathtools}
\usepackage{hyperref}

\newcommand{\R}{\mathbb{R}}
\newcommand{\Z}{\mathbb{Z}}

\newcommand{\be}{\begin{enumerate}}
\newcommand{\ee}{\end{enumerate}}
\newcommand{\ei}{\end{itemize}}

\newcommand{\ba}{\begin{array}}
\newcommand{\ea}{\end{array}}
\newcommand{\bmat}{\left[\begin{array}}
\newcommand{\emat}{\end{array}\right]}
\newcommand{\bt}{\begin{thm}}
\newcommand{\et}{\end{thm}}
\newcommand{\bp}{\begin{proof}}
\newcommand{\ep}{\end{proof}}
\newcommand{\bprop}{\begin{prop}}
\newcommand{\eprop}{\end{prop}}
\newcommand{\bl}{\begin{lemma}}
\newcommand{\el}{\end{lemma}}
\newcommand{\bc}{\begin{cor}}
\newcommand{\ec}{\end{cor}}
\newcommand{\bd}{\begin{defn}}
\newcommand{\ed}{\end{defn}}

\newcommand{\vep}{\varepsilon}

\newcommand{\tx}{\textrm}

\begin{document}

%%%%%%%%%%%%%%%%%%%%% Publisher's Area please ignore %%%%%%%%%%%%%%
\catchline{}{}{}{}{}
%%%%%%%%%%%%%%%%%%%%%%%%%%%%%%%%%%%%%%%%%%%%%%%%%%%%%%%%%%%%%%%%%%%

\title{A lower bound on critical points of the electric potential of a knot}

\author{Max Lipton}

\address{Cornell University \\
Malott Hall, 301 Tower Road \\
Ithaca, NY 14850 \\
ml2437@cornell.edu}

\maketitle

\begin{abstract}
\noindent Consider a knot $K$ in $S^3$ with charge uniformly distributed on it. From the standpoint of both physics and knot theory, it is natural to try to understand the critical points of the potential and their behavior. \\

\noindent  We show the number of critical points of the potential is at least $2t(K) + 2$, where $t(K)$ is the tunnel number, defined as the smallest number of arcs one must add to $K$ such that its complement is a handlebody. The result is proven using Morse theory and stable manifold theory.\\
\end{abstract}



\section{Introduction}

Our novel problem of interest is to analyze the zeros of the electric field around a charged, knotted wire fixed in place. Let $K \subseteq \R^3 \subseteq S^3$ be a smooth knot parametrized by the curve $r(t),$ $t \in [0,2\pi]$ with $r(0) = r(2\pi)$. We will take the convention that $S^3$ is the union of $\R^3$ and a single compactifying point at infinity. Suppose $K$ is endowed with a uniform charge distribution. With a choice of units, the electric potential between a point $k \in K$ and a point charge at $x$ at a distance $R$ from $k$ is proportional to $R^{-1}$. It therefore makes sense to define the electric potential $\Phi: S^3 - K \to \R$, on the complement of $K$ by the line integral 

\noindent
\begin{equation}
\Phi(x) = \int_{k \in K} \frac{dk}{|x - k|} = \int_0^{2\pi} \frac{|r'(t)|}{|x - r(t)|}dt, \tx{ } x \in \R^3 - K.
\label{this}
\end{equation}

We set $\Phi(\infty) = 0$ to ensure smoothness. By differentiating under the integral sign with respect to $x$, we can see the electric potential is smooth and harmonic. The \textit{electric field} is defined by $E = - \nabla \Phi$. We want to describe the critical points of the potential (equivalently, the zeros of the electric field) and their behaviors. These represent equilibrium points where a charged particle at rest will continue to experience no electric force from the charge distribution. Some conventions use the negative of the potential, so that the electric field points towards the knot, but it is more convenient for our purposes to work with a nonnegative potential function.

Define the knot invariant $cp(K)$ to be the smallest number of critical points of the electric potential among all parametrizations in the knot isotopy class $[K]$. All parametrizations have a critical point at infinity, which we include in the count. We will now assume $r$ is a parametrization which yields a critical set of minimal size. We obtain a lower bound for the number of critical points of the electric potential based on a well known topological invariant called the tunnel number $t(K)$. The tunnel number was originally introduced by Clark \cite{clark}, and remains an active topic of knot theory research \cite{baker}. 

We now come to the main result of this article. Proving this theorem uses Morse theory and stable manifold theory. 

\begin{theorem}
For all knots $K$, $cp(K) \geq 2t(K) + 2$.
\end{theorem}

\section{Preliminary Definitions and Lemmas}

The theorems of Morse theory require us to work on a compact manifold, so in the sequel we will define the knot complement of $K$ in $S^3$ to be the complement of an open tubular neighborhood of $K$ of sufficiently small radius. The theorems from Morse theory we invoke will still hold on this compact manifold with boundary because the potential function is proper, with gradient transversely intersecting the boundary. We will still denote our domain by $S^3 - K$.

Many of the following definitions and results are standard and can be found in Nicolaescu \cite{nico} and Burde \cite{burde}. A \textit{critical point} of a smooth real valued function $f$ on a manifold $M$ is a point $p$ such that the differential $df(p)$ is zero. The \textit{critical set} of $f$ is the set of critical points, and is denoted $\tx{Crit}(f)$. We say $f$ is \textit{Morse} if its critical points are nondegenerate, which means the Hessian matrix $H(f)$ of second partial derivatives is nonsingular. The $\textit{index}$ of $p \in \tx{Crit}(f)$ is the number of negative eigenvalues of $H(f)$, which is invariant under the choice of local coordinates. We denote the set of critical points of index $i$ by $\tx{Crit}_i(f)$. If $f$ is fixed, we denote the index of $p$ by $\lambda(p)$ and we denote the number of critical points of index $i$ by $m_i$. In the space of all smooth real valued maps, under a suitable function space topology, the set of Morse functions is dense. Therefore, we may assume the electric potential $\Phi$ is Morse by adding a perturbation if necessary.

We write $W^S(p)$ and $W^U(p)$ to denote the \textit{stable and unstable manifolds} of $p \in \tx{Crit}(f)$ respectively. Recall that the stable manifold (resp. unstable manifold) of $p$ is the set of all points which flow to $p$ along the gradient vector field $\nabla f$ as time tends to infinity (resp. negative infinity). If $f$ is Morse, the dimension of $W^S(p)$ is $\lambda(p)$ and the dimension of $W^U(p)$ is $\tx{dim }M - \lambda(p)$.

The \textit{Morse Reconstruction Theorem} states the domain of a Morse function on a compact manifold can be expressed as a cell complex by attaching closed discs with dimensions given by the indices of the critical points. The attaching maps are obtained by a process known as surgery, but we will not need to discuss the attaching maps in any further detail.

The \textit{Morse inequalities} state for a fixed Morse function $f$ on a manifold $M$, $m_i \geq b_i$, where $b_i = \tx{dim } H_i(M)$ is the $i$th Betti number of our domain manifold. We will need the stronger result which states 

\noindent
\begin{equation}
\sum\limits_{i = 0}^{\tx{dim } M} (-1)^i m_i = \sum\limits_{i = 0}^{\tx{dim } M} (-1)^i b_i = \chi(M),
\label{this}
\end{equation}

\noindent with $\chi(M)$ denoting the Euler characteristic of $M$. See Nicolaescu, Corollary 2.3.3 \cite{nico}.

Before turning to the proof of Theorem 1.1, we need to prove a few preliminary lemmas.

\begin{lemma}
For all knots $K$, $H_i(S^3 - K) = \Z$, for $i = 0,1$ and $H_i(S^3 - K) = 0$ for $i \geq 2$.
\end{lemma}

\begin{proof}

See Rolfsen, Proposition 3.A.3 \cite{rolfsen}.
\end{proof}

From the homology of the knot complement, along with $(2.1)$, we can deduce the following lemma.

\begin{lemma}
For all knots $K$, and with $\Phi$ defined above, $m_1 - m_2 = 1$.
\end{lemma}

\begin{proof}
Equation $(2.1)$ states the Euler characteristic equals the alternating sum of the $m_i$'s. In other words, $\chi (S^3 - K) = \sum\limits_{i = 0}^3 (-1)^im_i$. From Lemma 2.1, we can see that $\chi (S^3 - K) = 1 - 1 = 0$. Since $\Phi$ is harmonic, every critical point has index $1$ or $2$, save for the point at infinity, which has index $0$. We can conclude $m_0 - m_1 + m_2 - m_3 = 0$, or equivalently, $m_1 - m_2 = 1$ as desired. 
\end{proof}

\begin{remark}
As $m_2 \geq 0$, $m_1 = 1 + m_2 \geq 1$. That is, there is always a critical point of index $1$.
\end{remark}

The next set of definitions and results are standard in $3$-manifold topology, and further exposition can be found in \cite{schleimer} and \cite{scharlemann}. A \textit{handlebody} is a topological space homotopic to the three dimensional ball with solid handles attached (by  ``attaching handles" we mean there are copies of $D^2 \times [0,1]$ where the boundary discs $D^2 \times 0$ and $D^2 \times 1$ are embedded on the boundary of the three-ball). Given a knot $K \subseteq S^3$, the \textit{tunnel number} $t(K)$ is the least number of arcs we must add to $K$ such that the complement in $S^3$ is a handlebody. A collection of arcs with this property is known as a \textit{tunneling}.

In the proof of Theorem 1.1, we shall use the above definition of the tunnel number, but there is an equivalent definition that is more visually intuitive. A \textit{Heegaard splitting} of a three-manifold $M$ is an embedding of a closed, compact, and orientable surface $H$ such that the interior and exterior of $H$ in $M$ are both handlebodies. We say the genus of $H$ is the \textit{genus of the splitting}.

\begin{theorem}
Let $H$ be an unknotted embedding of a genus $g$ surface in $S^3$. That is, let $H$ be the boundary of a tubular neighborhood around a wedge sum of $g$ unknotted circles. Then $H$ defines a Heegaard splitting.
\end{theorem}

\begin{theorem}
Genus $g$ Heegaard splittings of $S^3$ are unique up to isotopy.
\end{theorem}

The previous result is known as \textit{Waldhausen's Theorem}. Clearly, a tunneling of a knot defines a Heegaard splitting. Therefore, we can view the tunnel number of $K$ as the least number of arcs we must add to $K$ such that it is isotopic to a wedge sum of unknotted circles. We can immediately deduce that tunnelings and the tunnel number always exist.

\begin{lemma} 
Every smooth knot $K$ has a tunnel number.
\end{lemma}

\begin{proof}
Take a diagram of $K$ with finitely many crossings. Over each crossing, introduce an arc connecting the top and bottom strands. Collapse each arc so that the top and bottom strands intersect, and project onto the diagram's plane so that we are left with a wedge sum of say, $g$ circles. By Theorem 2.4, it follows that the complement in $S^3$ is also a handlebody with $g$ handles.\\
\end{proof}

\begin{remark}
We just proved the tunnel number is bounded above by the crossing number, the least number of crossings needed in a knot diagram of $K$.
\end{remark}

As some elementary examples, the tunnel number of the unknot is zero, and the tunnel number of the trefoil is one. Indeed, the tunnel number need not be the crossing number. For example, torus knots have tunnel number $1$, yet can have arbitrarily high crossing number. See Clark \cite{clark}.

\section{Proof of Theorem 1.1}

We now come to the proof of our main result. To prove the result, we construct a tunneling with $m_2$ arcs. This proves $cp(K) \geq m_2 \geq t(K)$. Then, by applying Lemma 2.2, we get that $cp(K) = m_0 + m_1 + m_2 + m_3  \geq 1 + (t(K) + 1) + t(K) + 0 = 2t(K) + 2$ as desired.

\subsection{Construction of the tunneling}

We will apply the \textit{Morse Rearrangement Lemma} to allow us to make some additional convenient assumptions about $\Phi$ without losing generality. A proof can be found in Nicolaescu, Chapter 2.4 \cite{nico}. The theorem states we can find a smooth $\hat{\Phi}: S^3 - K \to \R$ satisfying the following properties:

\begin{itemize}
    \item $\Phi$ and $\hat{\Phi}$ share the same critical points, and each critical point has the same index.
    
    \item Suppose $p$ and $q$ are distinct critical points. Then $\hat{\Phi}(p) \neq \hat{\Phi}(q)$. If $\hat{\Phi}(p) < \hat{\Phi}(q)$, then $\lambda(p) \leq \lambda(q)$.
    
    \item Inside of a neighborhood of each critical point, the gradient flows for $\Phi$ and $\hat{\Phi}$ are identical. 
    
    \item If $\gamma$ is an integral curve of $\nabla \hat{\Phi}$, then $\Phi(\gamma(t))$ is strictly increasing in $t$. A vector field with this property is called $\textit{gradient-like with respect to } \Phi$. 
    
    \item For $p,q \in \tx{Crit}(\hat{\Phi})$, $W^S(p)$ and $W^U(q)$ intersect transversely. Morse functions with this property are called \textit{Morse-Smale}.
\end{itemize}

This theorem allows us to perturb the values of the critical points so we can reorder them ascending by index without affecting the topological data encoded by the original potential. At this point, we are not necessarily working with the physical potential whose formula is given in $(1.1)$, but for simplicity we will still refer to the perturbation as $\Phi$.

Our rearrangement restricts the limiting behavior of trajectories. Let $\gamma: (-\infty, \infty) \to S^3 - K$ be a trajectory for $\Phi$. When $t \to -\infty$, $\Phi(\gamma(t))$ strictly decreases, but it is bounded below by $0$, by assumption. Should $\lim\limits_{t \to -\infty} \Phi(\gamma(t)) = 0$, then $\lim\limits_{t \to -\infty} \gamma(t) = \infty$, the point of infinity on $S^3$, because it is the only point in the knot complement with zero potential. Otherwise, $\lim\limits_{t \to -\infty} \gamma(t)$ is a critical point. Similarly, we can deduce that for $t \to \infty$, either $\gamma(t)$ tends to a critical point or $K$, since $\Phi(\gamma(t))$ is strictly increasing.

Should both ends of $\gamma$ be critical points, then we know the index of the critical point at $t = -\infty$ is less than or equal to that of the critical point at $t = \infty$.

Consider the critical points of index $2$. For $p^2_1, \dots, p^2_{m_2} \in \tx{Crit}_2(\Phi)$, let $\Gamma_i = W^U(p^2_i)$. Notice each $\Gamma_i$ is a union of two trajectories leaving $p^2_i$, since the unstable manifold has dimension $1$. Since $\Phi$ strictly increases along trajectories, we have that the trajectories will either tend to $K$ or to another critical point of index $2$ as $t \to \infty$. However, should either end of $\Gamma_i$ tend to a critical point $q \in \tx{Crit}_2(\Phi)$, then the corresponding trajectory will be a submanifold of $W^S(q)$. However, $W^S(q)$ and $\Gamma_i$ are $2$ and $1$-submanifolds respectively, and for them to intersect transversely as per our assumption, their intersection cannot be more than $0$ dimensions. Therefore, we conclude both ends of $\Gamma_i$ eventually reach the tubular neighborhood of the knot. Our tunneling is only concerned with the arc outside of the tubular neighborhood, so we can assume $\Gamma_i$ is defined only on a compact interval. 

\begin{figure}[th]
\centerline{\includegraphics[width=2.2in]{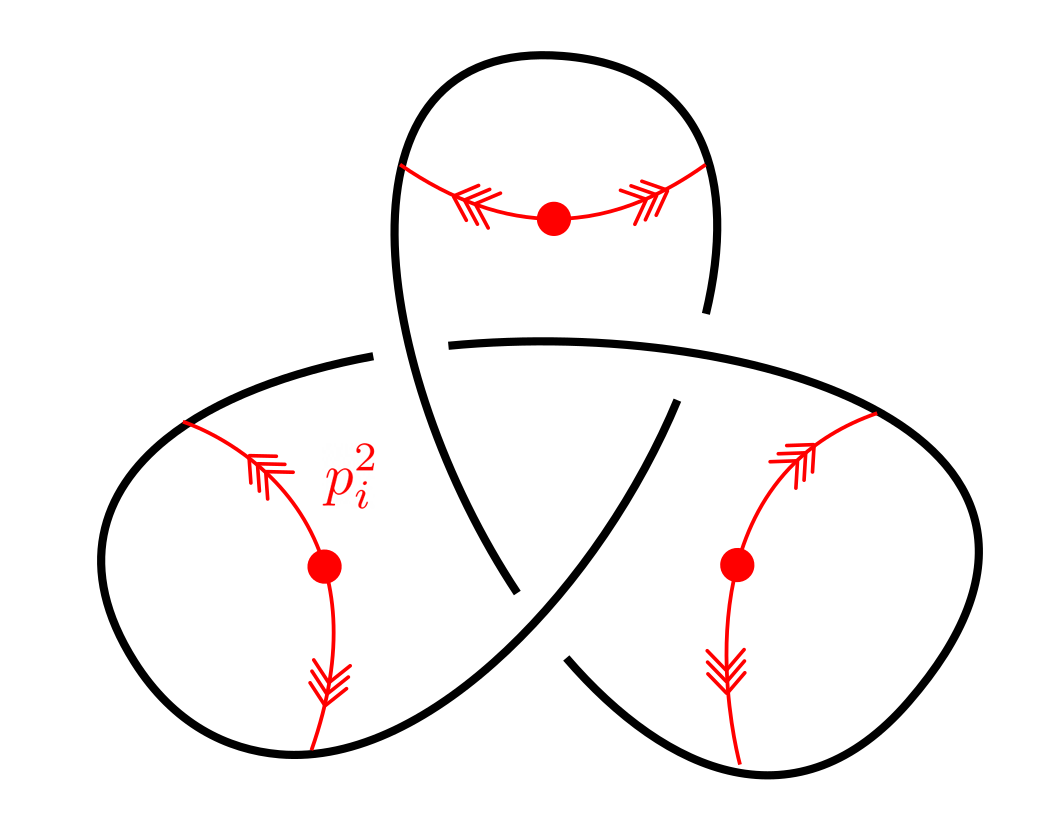}}
\vspace*{8pt}
\caption{The arcs we add to $K$ are the unstable manifolds associated to critical points of index $2$. Note that this diagram does not necessarily depict the specific situation accurately for the trefoil. \label{fig1}}
\end{figure}

Now consider the critical points of index $1$. For $p^1_1, \dots, p^1_{m_1} \in \tx{Crit}_1(\Phi)$, let $\Theta_j = W^S(p^1_j)$. Analogous to before, each $\Theta_j$ is a union of two trajectories tending to $p^1_j$. Similar reasoning will tell us that both ends tend to the point at infinity. Indeed, each $\Theta_j$ is a union of two trajectories tending towards a critical point $p^1_j$ of index $1$. As $t \to -\infty$, $\Phi$ will strictly decrease along these trajectories, so we know that the negative infinite limits of these trajectories must either be another critical point of index $1$, or the point at infinity. By the transversality assumption of the stable and unstable manifolds, we cannot have that the endpoints of $\Theta_j$ are critical points of index $1$. Therefore, both ends are at the point at infinity. We may view the union of all the $\Theta_j$ arcs as a wedge sum of circles at the point of infinity, which we denote $\bigvee\limits_{j = 1}^{m_1} \Theta_j$.

Notice $\bigvee\limits_{j = 1}^{m_1} \Theta_j$ is homotopy equivalent to a handlebody. Using a standard maneuver from differential topology, we will flow along the (negative) gradient to perform a deformation retraction from $S^3 - (K \cup \Gamma_1 \cup \dots \cup \Gamma_{m_2})$ to $\bigvee\limits_{j = 1}^{m_1} \Theta_j$. We will opt to work with smooth tubular neighborhoods of both of these spaces, but there is a crucial technical lemma we must prove before proceeding.

\subsection{Constructing a smooth boundary around the tunneling}

In this subsection, we prove the following lemma.

\begin{lemma}
There are tubular neighborhoods of both $K \cup \Gamma_1 \cup \dots \cup \Gamma_{m_2}$ and $\bigvee\limits_{j = 1}^{m_1} \Theta_j$ with smooth boundary such that the gradient vector field points inwards and outwards respectively.
\end{lemma}

\begin{proof}

\begin{figure}
\centerline{\includegraphics[width=4.5in]{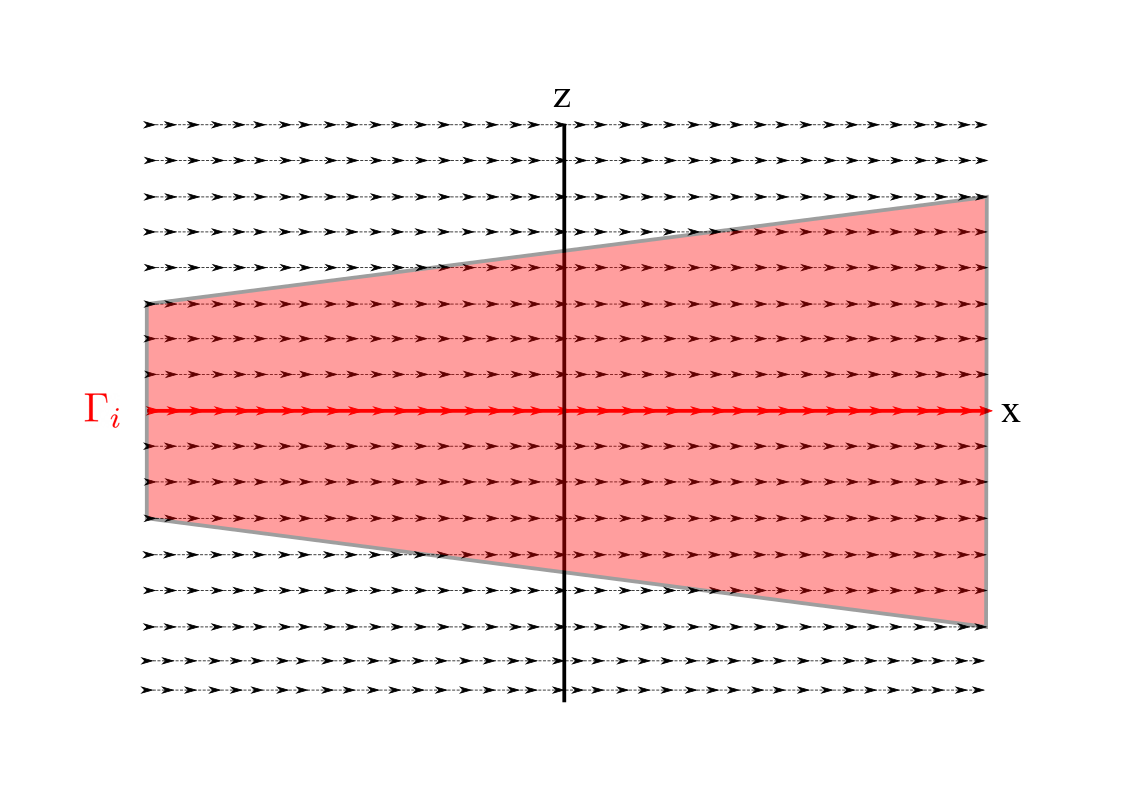}}
\caption{A sketch depicting the tubular neighborhood around a regular point of $\Gamma_i$ in specially chosen local coordinates which makes the gradient parallel. This is a projection to the $xz$-plane. \label{fig2}}
\end{figure}

\begin{figure}
\centerline{\includegraphics[width=3.5in]{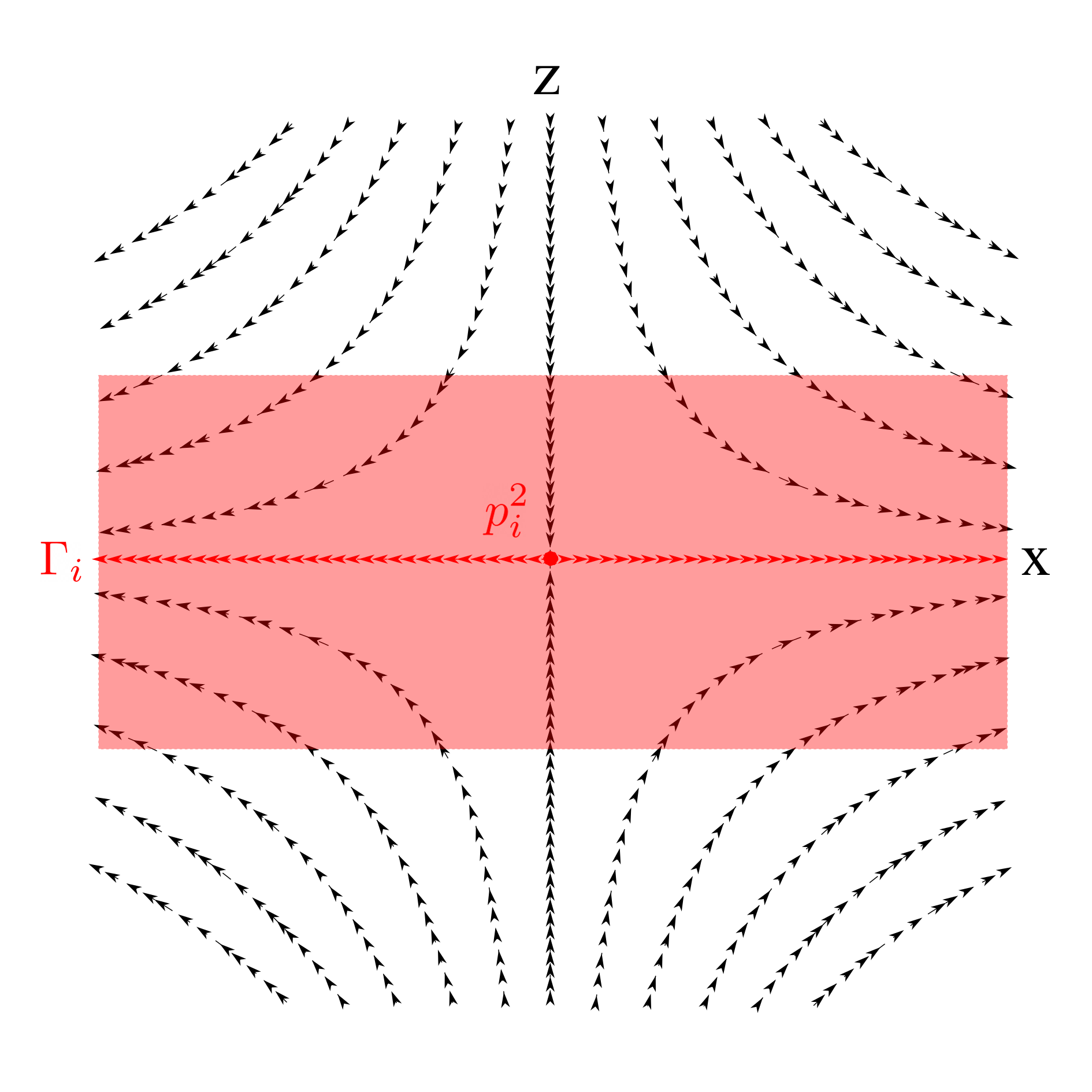}}
\caption{A sketch depicting the tubular neighborhood around a critical point on $\Gamma_i$, in specially chosen local coordinates such that $\Phi$ is a quadratic with signature $(1,2)$. This is a projection to the $xz$-plane. \label{fig3}}
\end{figure}

It is a standard fact from electrostatics that neighborhoods of $K$ and the point at infinity exist such that the gradient points inwards and outwards respectively. This proof will make use of the tubular flow lemma, which states for every regular point of $S^3 - K$, there are local coordinates such that the gradient flow takes the constant form $\nabla_{(x_0,y_0,z_0)} \Phi = (1,0,0) = \frac{\partial}{\partial x}$. Take local coordinates centered at a regular point of $\Gamma_i$ so that the portion of the arc inside our coordinate chart is mapped to the path of unit speed along the $x$-axis, $\gamma(t) = (t,0,0)$. 

Around a segment of the $x$-axis, we can choose our tubular neighborhood of $\Gamma_i$ to be the interior of the slanted tube defined by $y^2 + z^2 = (\frac{1}{2}x + 1)^2$, where we possibly restrict our local coordinates so $-1 < x < 1$. See Fig. 2. Notice the gradient points inwards from the boundary. If we expand around a regular point on a $\Theta_j$ arc, we can reflect the tube in the $x$ direction so the gradient points outwards. Note that the properties of a vector field pointing inwards and outwards from a boundary are invariant under a change of coordinate charts in an orientation-preserving atlas.

Now suppose we want to take local coordinates around a critical point $p^2_i$ of $\Gamma_i$, which has index $2$. By Lemma 2.2 of \cite{milnor}, we can find local coordinates centered at the critical point such that $\Gamma_i$ corresponds to the $x$-axis, and $\Phi$ takes the form $\Phi(x_0,y_0,z_0) = \frac{1}{2}(x_0^2 - y_0^2 - z_0^2) + c$, for some constant $c$. Thus, $\nabla \Phi$ takes the form $\nabla \Phi(x_0,y_0,z_0) = (x_0, -y_0, -z_0)$. Consider the tube around the $x$-axis defined by $y^2 + z^2 = 1$. See Fig. 3. Using an abuse of notation, we will refer to this tube as $\partial \Gamma_i$. We will show the gradient points inside the tube, as Figure 3 shows. At a point $(x_0,y_0,z_0) \in \partial \Gamma_i$, the tangent plane is spanned by $(1,0,0)$ and $(0,-z_0, y_0)$. As the point varies on the tube, this oriented basis of the tangent plane varies smoothly, thus defining an orientation of the tube. The triple of tangent vectors $\{(1,0,0), (0,-z_0, y_0), \nabla \Phi (x_0,y_0,z_0)\}$ therefore defines a smooth choice of orientation for the three dimensional ambient tangent spaces surrounding the tube. Therefore, if the gradient points inside the tube at one point of the tube, it does so throughout the whole tube. For example, at the point $p = (0,1,0) \in \partial \Gamma_i$, $\nabla \Phi(p) = (0,-1,0)$ clearly points inside the tube. The existence of a surface surrounding the arc with which the gradient points inwards is a significant topological obstruction to proving our main theorem. The construction would not be possible if the critical point had index $1$, or if we were asked to place a tube around another axis. The case for a critical point of $\Theta_i$ is analogous, and we get the result that the gradient points outwards.

By compactness of the arcs $\Gamma_i$ and $\Theta_j$ we only need to construct finitely many tubes around arc segments to encapsulate the entire arc. When the tubes cover the same part of the arc, we may have to shrink the radius of one of the tubes so the boundaries will intersect, but the resulting tube will still have the gradient pointing in the correct direction. By taking the final boundary to be the points of minimal radial distance to the arc, we obtain a connected boundary to the tubular neighborhood that is only piecewise smooth. Likewise, the tubular neighborhoods around the $\Gamma_i$ arcs intersect the tube around $K$, and the gradient points inwards on the boundary of the union. This construction is sufficient to prove the lemma. To get a smooth boundary out of the piecewise smooth boundary, one could either make a density argument in a space of manifolds \cite{nico}, justify why the theorems work in the piecewise smooth case \cite{milnor}, or use mollifiers to smooth out the kinks \cite{stein}.  For the sake of space, we omit the technical details. 

For reasons that will be clear when we perform the deformation retraction, we will want to include a ball around the point at infinity in the tubular neighborhood. In the standard $\R^3$ coordinates, this is the complement of a large open ball. We can assume this neighborhood around $\infty$ contains $\Phi^{-1}([0, \vep])$ for some $\vep > 0$. On the boundary sphere, we can again use mollifiers to connect the tube smoothly whilst preserving their orientations against the gradient flow.

\end{proof}

\subsection{The deformation retraction to a handlebody}

The final step is to use the flow of $E$ to perform the deformation retraction. We will use the closed tubular neighborhoods described in Lemma 3.1 as an alternative to just the knot with arcs attached and the wedge of circles. Let $A$ be the aforementioned tubular neighborhood of the knot with the $\Gamma_i$ arcs attached, and let $B$ denote the aforementioned tubular neighborhood of the $\Theta_j$ arcs connected to a ball around the point at infinity.

Let $\mathcal{E}(x,t)$ be the flow of the \textit{negative} gradient. The point $\mathcal{E}(x,t)$ refers to the location of the path at time $t$ starting from the unique integral curve starting at $x$. First, we prove every point in $S^3 - A$ will eventually flow to $B$. Suppose $x \in S^3 - A$. There are three possibilities for the limit of the negative gradient flow of $x$: it could flow to $K$, it could flow to a critical point of index $1$ or $2$, or it could flow to the point at infinity. Since the negative gradient points outwards from the boundary of $A$, $x$ cannot flow to $K$ or a critical point of index $2$. Therefore, $x$ either flows to a critical point of index $1$, or to the point at infinity, which means that $x$ eventually flows to $B$. Furthermore, since the negative gradient points into $B$, once $x$ enters $B$, it will never leave. This fact still holds even in the vacuous case where $m_2 = 0$ and therefore $A = K$.

By smoothness of the boundary of $B$, the function which assigns each $x \in S^3 - A$ the minimum time $t$ such that $\mathcal{E}(x,t) \in B$ is smooth. Call this function $C(x)$. Notice $C$ assigns $0$ to each point already in $B$. By compactness of the domain, the function reaches a maximum value $C_{max}$. Define a homotopy $H$ on $(S^3 - A) \times [0, \infty)$ by $H(x,t) = \mathcal{E}(x, \tx{min}(t,C(x)))$, which we can see is continuous. Also notice that for $x \in B$, $H(x,t) = x$ for all $t$. Running this homotopy on the time interval $[0,C_{max}]$ completes the deformation retraction. This completes the proof of our main theorem. $\square$

\section*{Acknowledgments}
The author would like to thank his advisor, Steven Strogatz, for introducing him to this problem, and for helping motivate the connections to physics. He would also like to thank Greg Buck, Rennie Mirollo, Tim Riley, Alex Townsend, and the journal reviewer for their helpful discussions, comments, and suggestions. Finally, he would like to thank the participants of the RTG Math Communications Seminar (of whom there are too many to list) for their feedback on the writing.

This work was funded in part by the NSF RTG Grant entitled Dynamics, Probability, and Partial Differential Equations in Pure and Applied Mathematics, DMS-1645643.

\end{document}